 \documentclass[11pt]{article}
 \usepackage[top=1in,bottom=1in,left=1.5in,right=1.5in]{geometry}
  \usepackage{amsmath,amssymb,amscd}
 \usepackage{latexsym}
 \usepackage[dvips]{pstcol} 
  \usepackage{pst-node}
  \usepackage[dvips]{graphicx}

  \renewcommand{\epsilon}{\varepsilon}
  \newcommand{\C}{\mathbb{C}}
  
  \newcommand{\N}{\mathbb{N}}
  \renewcommand{\P}{\mathbb{P}}
  \newcommand{\R}{\mathbb{R}}
  \newcommand{\Z}{\mathbb{Z}}

  \renewcommand{\a}{\mathbf{a}}
  \renewcommand{\b}{\mathbf{b}}

  \newcommand{\w}{\mathbf{w}}
  \newcommand{\x}{\mathbf{x}}
  
  \newcommand{\y}{\mathbf{y}}
  \newcommand{\z}{\mathbf{z}}
  \newcommand{\0}{\mathbf{0}}

  \newcommand{\cB}{\mathcal{B}}
  
  \newcommand{\cD}{\mathcal{D}}

  \newcommand{\cX}{\mathcal{X}}
  \newcommand{\cY}{\mathcal{Y}}

  \newcommand{\rK}{\mathrm{K}}

  \newcommand{\lan}{\langle}
  \newcommand{\ran}{\rangle}
  \newcommand{\an}[1]{\lan#1\ran}

  \newcommand{\orb}{\mathrm{orb\;}}

  \newcommand{\rdyn}{\rho_{\mathrm{dyn}}}

  \newcommand{\Cl}{\mathrm{Cl\;}}
  
  \newcommand{\dist}{\mathrm{dist}}

  \newcommand{\Sing}{\mathrm{Sing\;}}

  \newcommand{\rC}{\mathrm{C}}
  \newcommand{\rd}{\mathrm{d}}
  \newcommand{\rH}{\mathrm{H}}

  \newcommand{\rS}{\mathrm{S}}

  \newcommand{\vol}{\text{vol}}
  \newcommand{\lov}{\text{lov}}

  \newtheorem{theo}{Theorem}[section]

  \newtheorem{con}[theo]{Conjecture}

  \numberwithin{equation}{section} 

\begin{document}

 \title{Entropy of holomorphic and rational maps: a survey}

 \author
 {Shmuel Friedland\\
 Department of Mathematics, Statistics and
 Computer Science\\
 University of Illinois at Chicago\\
 Chicago, Illinois 60607-7045
 }
 \date{July 24, 2006}
 \maketitle

 \begin{abstract}
 We give a brief survey on the entropy of holomorphic self maps $f$
 of compact K\"ahler manifolds, and rational dominating self maps
 $f$ of smooth projective varieties.  We emphasize the connection
 between the entropy and the spectral radii of the induced action of
 $f$ on the homology of the compact manifold.  The main
 conjecture for the rational maps states that modulo birational
 isomorphism all various notions of
 entropy and the spectral radii are equal.
     \\[\baselineskip] 2000 Mathematics Subject
     Classification: 28D20, 30D05, 37F, 54H20.

 \end{abstract}

 \section{Introduction}

 The subject of the dynamics of a map $f: X\to X$ has been
 studied by hundreds, or perhaps thousands, of mathematicians, physicists and
 other scientists in the last 150 years.  One way to classify the \emph{complexity}
 of the map $f$ is to assign to it a number $h(f)\in [0,\infty]$,
 which called the \emph{entropy} of $f$.  The entropy of $f$ should
 be an invariant with respect to certain \emph{automorphisms} of
 $X$.  The complexity of the dynamics of $f$ should be reflected by $h(f)$,
 i.e. the larger $h(f)$ the more complex is its dynamics.

 The subject of this short survey paper is mostly concerned with the entropy of
 a holomorphic $f:X\to X$, where $X$ is a compact K\"ahler
 manifold,  and the entropy of a rational map of $f: Y \dashrightarrow Y$, where
 $Y$ is a smooth projective variety.  In the holomorphic case the
 author \cite{Fr1,Fr4,Fr2} showed that entropy of $f$ is equal to
 the logarithm of the spectral radius of the finite dimensional
 $f_*$ on the total homology group $\rH_*(X)$ over $\R$.

 Most of the paper is devoted to the rational map $f:Y\dashrightarrow
 Y$ which can be assumed dominating.  In this case we have some
 partial results and inequalities.  We recall three
 possible definition of the entropy $h_B(f),h(f),h_F(f)$ which are
 related as follows: $h_B(f)\le h(f)= h_F(f)$.  The analog of the
 \emph{dynamical} homological spectral radius are given by
 $\rdyn(f_*)$, $e^{ \lov (f)}$ and $e^{H(f)}$, where the three quantities
 can be viewed as the volume growth.
 It is known that $h_F(f)\le \lov(f)\le H(f)$.
 $H(f)$ is a birational invariant.  I.e. let $ \hat Y$ be a smooth
 projective variety such that there exists $\iota:Y \dashrightarrow \hat Y$
 which is a \emph{birational} map.  Then $f:Y\dashrightarrow Y$ can be lifted
 to dominating $\hat f:=\iota f\iota^{-1}:\hat Y\dashrightarrow \hat Y$, and $H(f)=H(\hat f)$.
 However $h_F(f)$ does not have to be equal to $h_F(\hat f)$.
 The main conjecture of this paper are the equalities
 \begin{equation}\label{entropconeq}
 h_B(\hat f)=h(\hat f)=h_F(\hat f)=\lov(\hat f)=H(\hat f)=\rdyn(\hat f_{*}),
 \end{equation}
 for some $\hat f$ birationally equivalent to $f$.
 For polynomial automorphisms of $\C^2$, which are birational maps
 of $\P^2$, the results of the papers \cite{FM, Smi, DS} prove the above
 conjecture for $\hat f=f$.  Some other examples where this conjecture
 holds are given in \cite{Gue1,Gue2}.

 The pioneering inequality of Gromov $h_F(f)\le \lov(f)$ \cite{Gro} uses
 basic results in entropy theory, Riemannian geometry and complex manifolds.
 Author's results are using basic results in entropy theory,
 algebraic geometry and the results of Gromov, Yomdin \cite{Yom} and
 Newhouse \cite{New}.  From the beginning of 90's the notion of
 \emph{currents} were introduced in the study of the dynamics of
 holomorphic and rational maps in several complex variables.
 See the survey paper \cite{Sib}.  In fact the inequality
 $\lov(f)\le H(f)$ proved in \cite{DS,DS1,DS2} and \cite{Gue},
 as well as most of the
 results in  are derived \cite{Gue1,Gue2}, are
 using the theory of currents.

 The author believes that in dealing with the notion of the
 entropy solely, one can cleverly substitute the theory of currents with
 the right notions of algebraic geometry.  All the section of this
 paper except the last one are not using currents.
 It seems to the author that to prove the conjecture
 (\ref{entropconeq}) one needs to prove a correct analog of Yomdin's
 inequality \cite{Yom}.

 We now survey briefly the contents of this paper.
 \S2 deals with the entropy of $f:X\to X$, where first $X$
 is a compact metric space and $f$ is continuous, and second
 $X$ is compact K\"ahler and $f$ is holomorphic.  \S3 is devoted to
 the study of three definitions of entropy of a continuous map
 $f:X\to X$, where $X$ is an arbitrary subset of a compact metric
 space $Y$.  In \S4 we discuss
 rational dominating maps $f:Y\dashrightarrow Y$,
 where $Y$ is a smooth projective variety.
 \S5 discusses various notions and results on the entropy of
 rational dominating maps.
 In \S6 we discuss briefly the recent results, in particular the inequality $\lov(f)\le H(f)$
 which uses currents.

 It is impossible to mention all the relevant existing literature,
 and I apologize to the authors whose papers were not mentioned.
 It is my pleasure to thank S. Cantat, V. Guedj, J. Propp, N. Sibony
 and C.-M. Viallet for
 pointing out related papers.

 \section{Entropy of continuous and holomorphic maps}

 The first rigorous definition of the entropy was introduced
 by Kolmogorov \cite{Kol}.  It assumes that $X$ is a probability
 space $(X,\cB,\mu)$, where $f$ preserves the probability measure
 $\mu$.  It is denoted by $h_{\mu}(f)$, and is usually referred under
 the following names: \emph{metric} entropy,
 \emph{Kolmogorov-Sinai} entropy, or \emph{measure entropy}.
 $h_{\mu}(f)$ is an invariant under measure preserving invertible
 automorphism $A: X\to X$, i.e. $h_{\mu}(f)=h_{\mu}(A\circ f\circ A^{-1})$.

 Assume that $X$ is a compact metric space and $f:X\to X$ a continuous map.
 Then Adler, Konheim and McAndrew defined the \emph{topological} entropy $h(f)$
 \cite{AKM}.  $h(f)$ has a \emph{maximal characterization} in terms of
 measure entropies $f$.  Let $\cB$ be the Borel sigma algebra
 generated by open set in $X$.  Denote by $\Pi(X)$ the compact
 space of probability measures on $(X,\cB)$.  Let $\Pi(f)\subseteq
 \Pi(X)$ be the compact set of all $f$-invariant probability measures.
 (Krylov-Bogolyubov theorem implies that $\Pi(f)\ne \emptyset$.)
 Then the \emph{variational principle} due to
 Goodwyn, Dinaburg and Goodman \cite{Goo,Din,Gdm} states
 $h(f)=\max_{\mu\in\Pi(f)} h_{\mu}(f)$.  Hence $h(f)$ depends
 only to the topology induced by the metric on $X$.
 In particular, $h(f)$ is invariant under any homeomorphism $A:X\to
 X$.

 The next step is to consider the case where $X$ is a compact smooth
 manifold and $f:X\to X$ is a differentiable map, i.e. $f\in
 \rC^{r}(X)$, where $r$ is usually at least $1$.  The most remarkable
 subclasses of $f$ are strongly hyperbolic maps, and in particular
 axiom A diffeomorphisms \cite{Shu}.  The dynamics of an Axiom A diffeomorphism
 on the nonwandering set
 can be coded as a \emph{subshift of a finite type} (SOFT), hence
 its entropy is given by the exponential growth of the periodic
 points of $f$, i.e. $h(f)=\limsup_{k\to\infty} \frac{\log \textrm{Fix }f^k}{k}$,
 where  $\textrm{Fix }f^k$
 the number of periodic points of $f$ of period $k$.

 It is well known in topology that
 $\textrm{Fix }f^k$ can be estimated below by
 the Lefschetz number of $f^k$.  Let $\rH_*(X)$ denote
 the total homology group of $X$ over $\R$, i.e.
 $\rH_*(X)=\oplus_{i=0}^{\dim_{\R} X} \rH_i(X)$,
 the direct sum of
 the homology groups of $X$ of all dimensions with coefficients in $\R$.
 Then $f$ induces the linear operator $f_*: \rH(X) \to
 \rH(X)$, where $f_{*,i}:\rH_i(X)\to \rH_i(X), i=0,\ldots, \dim_{\R}
 X$.  The \emph{Lefschetz number} of $f^k$ is defined as $\Lambda(f^k):=\sum_{i=0}^{\dim_{\R}
 X} (-1)^i \textrm{Trace} f^k_{*,i}$.  Intuitively, $\Lambda(f^k)$
 is the algebraic sum of $k$-periodic points of $f$, counted with their multiplicities.

 Denote by $\rho(f_*)$ and $\rho(f_{*,i})$ the spectral radius
 of $f_*$ and $f_{*,i}$ respectively.
 Recall that $\rho(f_{*,i})=\limsup_{k\to\infty} |\textrm{Trace}
 f^k_{*,i}|^{\frac{1}{k}}$ and
 $\rho(f_*)=\max_{i=0,\ldots,\dim_{\R} X} \rho(f_{*,i})$.
 Hence
 $$\limsup_{k\to\infty} \frac{\log |\Lambda(f^k)|}{k}\le
 \log\rho(f_*).$$
 The arguments in \cite{Shu} yield that for any $f$ in the subset $H$
 of an Axiom A diffeomorphism, ($H$ is defined in \cite{Shu}), one has the inequality
 $|\textrm{Trace}
 f^k_{*,i}|\le \textrm{Fix }f^k$ for each $=1,\ldots, \dim_{\R} X$.  ($H$ is $C^0$  dense in
 $\textrm{Diff}^r(X)$ \cite[Thm 3.1]{Shu}.)
 Hence for any $f\in H$ one has the inequality \cite[Prop 3.3]{Shu}
 \begin{equation}\label{shubcon}
 h(f)\ge \log\rho(f_*).
 \end{equation}
 It was conjectured in \cite{Shu} that the above inequality holds
 for any differentiable $f$.

 Let $\deg f$ be the \emph{topological degree} of $f:X\to X$.
 Then $|\deg f|=\rho(f_{*,\dim_{\R} X})$.
 Hence $\rho(f_*)\ge |\deg
 f|$.  It was shown by Misiurewicz and Przytycki \cite{MP} that if
 $f\in \rC^{1}(X)$ then $h(f)\ge |\deg f|$.  However this
 inequality may fail if $f\in\rC^0(X)$.
 The entropy conjecture (\ref{shubcon}) for  a smooth $f$, i.e. $f\in \rC^{\infty}(X)$,
 was proved by Yomdin \cite{Yom}.  Conversely, Newhouse \cite{New} showed that
 for $f\in \rC^{1+\epsilon}(X)$ the \emph{volume growth} of smooth submanifolds of $f$
 is an upper bound for $h(f)$.  See also a related upper bound in
 \cite{Prz}.

 This paper is devoted to study the entropy of $f$ where $X$ is a
 complex K\"ahler manifold and $f$ is either holomorphic map, or $X$
 is a projective variety and $f$ is a rational map dominating map.
 We first discuss the case where $f$ is holomorphic.

 Let $\P$ be the complex projective space.  Then $f:\P\to\P$ is
 holomorphic if and only if $f|\C$ is a rational map.  Hence $\deg
 f$ is the cardinality of the set $f^{-1}(z)$ for all but a finite
 number of $z\in\C$.  So $\deg f=\rho(f_*)$ in this case.
 Lyubich \cite{Lyu} showed that $h(f)=\log\deg f$.
 Gromov in preprint dated 1977, which  appeared as \cite{Gro},
 showed that if $f:\P^d\to \P^d$ is holomorphic then
 $h(f)=\log\deg f$.  It is well known in this case
 $\rho(f_*)=\deg f$.

 In \cite{Fr1} the author showed that if $X$ is a complex projective
 variety and $f:X\to X$ is holomorphic, then $h(f)=\log\rho(f_*)$.
 Note that one can view $f_*$ a linear operator on $H_*(X,\Z)$,
 i.e. the total homology group with integer coefficients.
 Hence $f_*$ can be represented by matrix with integer coefficients.
 In particular, $\rho(f_*)$ is an algebraic integer, i.e. the entropy
 is the logarithm of an algebraic integer.
 (This fact was observed in \cite{BV} for certain rational maps.)
 In \cite{Fr2} the author extended this result to a
 compact K\"ahler manifold.

 Examples of the dynamics of biholomorphic maps
 $f:X\to X$, where $X$ is a compact $K3$ surface
 which is K\"ahler but not necessary a
 projective variety, are given in \cite{Can,McM}.  See also \cite{DS1}
 for higher dimensional examples.
 In summary, the entropy of a holomorphic self map $f$ of a compact
 K\"ahler manifold is determined by the spectral radius of
 the induced action of $f$ on the total homology of $X$.

 \section{Definitions of entropy}

 In this paper $Y$ will be always a compact matrix space with the metric
 $\dist(\cdot,\cdot): Y\times Y\to \R_+$.  Let $X\subseteq Y$ be a
 nonempty set, and assume that $f:X\to X$ is a continuous map
 with respect the topology induced by the metric $\dist$ on $X$.
 For $x,y\in X$ and $n\in\N$ let
 $$\dist_n(x,y)=\max_{k=0,\ldots,n-1} \dist(f^k(x),f^k(y)).$$
 So $\dist_1(x,y)=\dist(x,y)$ and the sequence $\dist_n(x,y), n\in\N$
 is nondecreasing.  Hence for each $n\in\N$ $\dist_n$ is a distance
 on $X$.  Furthermore, each metric $\dist_n$ induces the same topology $X$ as
 the metric $\dist$.  For $\epsilon>0$ a set $S\subseteq X$ is called
 $(n,\epsilon)$ separated if $\dist_n(x,y)\ge \epsilon$ for any $x,y\in S, x\ne
 y$.  For any set $K\subseteq X$ denote by $N(n,\epsilon,K)\in \N\cup\{\infty\}$ the maximal
 cardinality of $(n,\epsilon)$ separated set $S\subseteq K$.  Clearly,
 $N(n_1,\epsilon,K)\ge N(n_2,\epsilon,K)$ if $n_1 \ge n_2$,
 $N(n,\epsilon_1,K)\ge
 N(n,\epsilon_2,K)$ if $0<\epsilon_1\le \epsilon_2$, and $N(n,\epsilon,K_1)\ge
 N(n,\epsilon,K_2)$ if $K_1\supseteq K_2$.

 We now discuss a few possible definitions of the entropy of $f$.
 Let $K\subseteq X$.  Then
 \begin{equation}\label{deftopencx}
 h(f,K):=\lim_{\epsilon \searrow 0} \limsup_{n\to\infty} \frac{\log
 N(n,\epsilon,K)}{n}.
 \end{equation}
 We call $h(f,K)$ the \emph{topological entropy} of $f|K$.
 (Note that $h(f,K)=\infty$ if $N(n,\epsilon,K)=\infty$ for some
 $n\in\N$ and $\epsilon>0$.)  Equivalently, $h(f,K)$ can be viewed
 as the \emph{exponential growth} of the maximal number of
 $(n,\epsilon)$ separated sets (in $n$).

 Clearly $h(f,K_1)\ge h(f,K_2)$ if $X\supseteq K_1\supseteq K_2$.
 Then $h(f):=h(f,X)$ is the \emph{topological entropy} of $f$.

 Bowen's definition of
 the entropy of $f$, denoted here as $h_B(f)$, is given as follows \cite[\S7.2]{Wal}.
 Let $K\subseteq X$ be a compact set.  Then $K$ is a compact set with
 respect to $\dist_n$.  Hence $N(n,\epsilon,K)\in\N$.
 Then $h_B(f,X)$ is the supremum of $h(f,K)$ for all
 compact subsets $K$ of $X$.
 I.e.
 $$h_B(f,X)=\sup_{K\Subset X} h(f,K).$$
 When no ambiguity arises we let
 $h_B(f):=h_B(f,X)$.
 Clearly, if $X$ is compact then $h_B(f)=h(f)$.
 (It is known that for a compact $X$ $h(f)\in[0,\infty]$, i.e.
 \cite{Wal}.)

 Since $N(n,\epsilon,K)\le N(n,\epsilon,X)$ for any $K\subseteq X$
 it follows that $h(f)\ge h_B(f)$.  The following example, pointed out to me
 by Jim Propp, shows that
 it is possible that $h(f)>h_B(f)$.  Let $Y:=\{z\in\C,\;|z|\le 1\}, X:=
 \{z\in\C,\;|z|< 1\}$ be
 the closed and the open unit disk respectively in the complex plane.
 Let $2\le p\in \N$ and assume that $f(z):=z^p$.  It is well known
 that $h(f,Y)=\log p$.  It is straightforward to show that
 $h(f,X)=h(f,Y)$.  Let $K\subset X$ be a compact set.  Let $D(0,r)$ be
 the closed disk or radius $r<1$, centered at $0$, such that $K\subseteq
 D(0,r)$. Since $f(D(0,r))\subseteq D(0,r)$ it follows that $h_B(f,X)\le h(f,D(0,r))=0$.

 Our last definition of the entropy of $h$, denoted by $h_F(f,X)$,
 or simply $h_F(f)$ is based on the notion of the orbit space.
 Let $\cY:=Y^{\N}$ be the space of the sequences
 $\y=(y_i)_{i\in\N}$, where each $y_i\in Y$.  We introduce a metric
 on $\cY$:
 $$\rd(\{x_i\},\{y_i\}):=\sum_{i=1}^{\infty}
 \frac{\dist(x_i,y_i)}{2^{i-1}}, \quad \{x_i\}_{i\in\N},
 \{y_i\}_{i\in\N}\in \cY.$$
 Then $\cY$ is a compact metric space, whose diameter is twice the diameter of $Y$.
 The \emph{shift} transformation $\sigma: \cY\to \cY$ is given by
 $\sigma(\{y_i\}_{i\in\N})=\{y_{i+1}\}_{i\in\N}$.
 Then
 $\rd(\sigma(\x),\sigma(\y))\le 2 \rd(\x,\y)$, i.e.
 $\sigma$ is a Lipschitz map.
 Given $x\in X$ then the $f$-\emph{orbit} of $x$, or simply the orbit of $x$,
 is the point
 $\orb x:=\{f^{i-1}(x)\}_{i\in \N}\in\cY$.  Denote by $\orb X\subseteq \cY$,
 the \emph{ orbit space}, the
 set of all $f$-orbits.  Note that $\sigma(\orb x)=\orb f(x)$.
 Hence $\sigma(\orb X) \subseteq \orb X $.  $\sigma |\orb X$, the
 restriction of $\sigma$ to the orbit space, is ``equivalent" to the
 map $f: X\to X$.  I.e. let $\omega:X \to \orb X$ be given by
 $\omega(x):=\orb x$.  Clearly $\omega$ is a homeomorphism.
 Then the following diagram is commutative:
 \[
 \begin{CD} X  @>f>> X\\ @V\omega VV @ V\omega VV\\
 \cX @>\sigma>> \cX\\
 \end{CD}
 \]

 Let $\cX$ be the closure of $\orb X$ with respect
 to the metric $\rd$ defined above.  Since  $\cY$ is compact, $\cX$
 is compact.  Clearly $\sigma (\cX)\subseteq \cX$.  Following
 \cite[\S4]{Fr1} we define $h_F(f,X)$ to be equal to the topological
 entropy of $\sigma|\cX$:
 $$h_F(f,X):=h(\sigma|\cX)=h(\sigma,\cX).$$
 When no ambiguity arises we let
 $h_F(f):=h_F(f,X)$.
 Since the closure of $\orb X$ is $\cX$, it is not difficult to show
 that $h_F(f)=h(\sigma,\orb X)$.

 Observe first that if $X$ is a compact subset of $Y$ then
 $h_F(f)$ is the topological entropy $h(f)$ of $f$.  Indeed,
 since $f$ is continuous and $X$ is compact $\cX=\orb X$.
 Since $\omega$ is a homeomorphism,
 the variational principle implies that $h(f)=h_F(f)$.

 We observe next that $h(f)\le h_F(f)$.
 Let
 $$\rd_n(\x,\y):=\max_{k=0,\ldots,n-1}
 \rd(\sigma^{k}(\x),\sigma^{k}(\y)).$$
 Then $\dist_n(x,y) \le \rd_n(\orb x,\orb y)$.
 Hence $N(n,\epsilon,X)\le N(n,\epsilon,\cX)$.
 Hence $h(f)\le h_F(f)$.  The arguments of the proof \cite[Lemma
 1.1]{Gue2} show that $h(f)=h_F(f)$.  (In \cite{Gue2}
$h^{\textrm{Bow}}_{\textrm{top}}(f)$
 is our $h(f)$, and $h^{\textrm{Gr}}_{\textrm{top}}(f)$ is the
 topological entropy with respect to the metric
 $\rd'(\{x_i\},\{y_i\}):=\sup_{i\in\N}
 \frac{\dist(x_i,y_i)}{2^{i}}$.  Since $\rd$ and $\rd'$ induce the
 Tychonoff topology on $Y^{\N}$ it follows that
 $h^{\textrm{Gr}}_{\textrm{top}}(f)=h_F(f)$.)

 Our discussion of various topological entropies for $f:X\to X$
 is very close to the discussion in \cite{HNP}.
 The notion of the entropy $h_F(f)$ can be naturally extended to
 the definition of the entropy of a semigroup acting on $X$
 \cite{Fr3}.  See \cite{Buf} for other definition of the entropy
 of a free semigroup and \cite{Fis} for an analog of
 Misiurewicz-Przytycki theorem \cite{MP}.

 \section{Rational maps}

 In this section we use notions and results from algebraic geometry
 most of which can be found in \cite{GH}.
 Let $\z=(z_0,z_1,\ldots,z_n)$, sometimes denotes as $(z_0:z_1:\ldots:z_n)$,
 be the homogeneous coordinates the $n$-dimensional complex projective space $\P^n$.
 Recall that a map $f:\P^n \dashrightarrow \P^n$ is called a rational map if
 there exists $n+1$ nonzero coprime homogeneous polynomials
 $f_0(\z),\ldots,f_n(\z)$ of degree $d\in\N$ such that $\z\mapsto
 f_h(\z):=(f_0(\z),\ldots,f_n(\z))$.
 Equivalently $f$ lifts to a homogeneous map
 $f_h:\C^{n+1}\to\C^{n+1} $.
 The set of singular points of
 $f$, denoted by $\Sing f \subset \P^n$,
 sometimes called the \emph{indeterminacy locus} of
 $f$, is given by the system $f_0(\z)=\ldots=f_n(\z)=0$.
 $\Sing f$ is closed subvariety of $\P^n$ of codimension $2$
 at least.  $f$ is holomorphic if and only if $\Sing f=\emptyset$,
 i.e. the above system
 of polynomial equations has only the solution $\z=\0$.

 Let $Y$ be an irreducible algebraic variety.
 It is well known that $Y$ can be embedded as
 an irreducible subvariety of $\P^n$.  For simplicity of notation
 we will assume that $Y$ is an irreducible variety of $\P^n$.
 So $Y$ can be viewed as a homogeneous irreducible variety $ Y_h\subset
 \C^{n+1}$, given as the zero set of homogeneous polynomials
 $p_1(\z)=\ldots=p_m(\z)=0$.  $y\in Y$ is called \emph{smooth} if
 $Y$ is a complex compact manifold in the neighborhood of $y$.
 A nonsmooth $y\in Y$ is called a \emph{singular} point.
 The set of singular points of $Y$, denoted by $\Sing Y$,
 is a strict subvariety of $Y$.  $Y$ is called \emph{smooth} if
 $\Sing Y=\emptyset$.  Otherwise $Y$ is called \emph{singular}.

 Let $f:Y \dashrightarrow Y$ be a
 rational map.
 Then one can extend $f$
 to $\underline{f} :\P^n \dashrightarrow \P^n$ such that $\Sing
 \underline {f}
 \cap Y$ is a strict subvariety of $Y$ and
 $\underline{ f}|(Y\backslash \Sing \underline{ f})=
 f|(Y\backslash \Sing \underline {f})$.  $\underline{ f}$ is not unique, but the
 $f$ can be viewed as $\underline{ f}|Y$.  $\Sing f \subset Y$ is the
 set of the points where $f$ is not holomorphic.  $\Sing f$ is
 strict projective variety of $X$, ($\Sing f\subseteq \Sing \underline{
 f} \cap Y$), and each irreducible component of $\Sing f$ is at least of codimension
 $2$.  The assumption $f:Y
 \dashrightarrow Y$ means that $\w:=f_h(\z)\in Y_h$ for each
 $\z\in Y_h$.  It is known that $Y_1:=\Cl f_h(Y_h)$, the closure
 of $f_h(Y_h)$, is a homogeneous irreducible subvariety of $Y$.
 Furthermore either $Y_1=Y(=Y_0)$, in this case $f$ is called a
 \emph{dominating} map, or $\dim Y_1 < \dim Y_0$.  In the second case
 the dynamics of $f_0:=f$ is reduced to the dynamics of the rational map
 $f_1:Y_1\dashrightarrow Y_1$.  Continuing in the same manner we
 deduce that there exists a finite number of strictly descending
 irreducible subvarieties $Y_0:=Y \supsetneqq \ldots \supsetneqq Y_k$
 such that $f_k:Y_k\dashrightarrow Y_k$ is a rational dominating
 map.  (Note that $Y_k$ may be a singular variety.)
 Thus one needs only to study the dynamics of a rational
 dominating map $f:Y \dashrightarrow Y$, where $Y$ may be a singular variety.

 The next notion is the resolution of singularities of $Y$ and $f$.
 An irreducible projective variety $Z$ birationally equivalent to
 $Y$ if the exists a birational map $\iota:Z \dashrightarrow Y$.
 $Z$ is called a \emph{blow up} of $Y$ if there exists a birational
 map $\pi:Z \to Y$ such $\pi$ is holomorphic.   $Y$ is called a blow
 down of $Z$.  Hironaka's result claims that  any irreducible
 singular variety $Y$ has a smooth blow up $Z$.  Let
 $f:Y\dashrightarrow Y$ be a rational dominating map.  Let $Y$ be a
 birationally equivalent to $Z$.  Then $f$ lifts to a rationally
 dominating map $\hat f:Z \dashrightarrow Z$.  Hence to study the dynamics
 of $f$ one can assume that $f:Y \dashrightarrow Y$ is rational
 dominating map and $Y$ is smooth.
 Hironaka's theorem implies that there exists a smooth blow up $Z$ of
 $Y$ such that $f$ lifts to a holomorphic map $\tilde f:Z \to Y$.
 Then one has the induced dual linear maps on the homologies and the
 cohomologies of $Y$ and $Z$:
 $$\tilde f_*:\rH_*(Z) \to \rH_*(Y), \quad \tilde f^*: \rH^*(Y) \to
 \rH^*(Z).$$
 We will view the homologies $\rH_*(Y), \rH_*(Z)$ as homologies with coefficients in
 $\R$, and hence the cohomologies
 $\rH^*(Y),\rH^*(Z)$, which are dual to $\rH_*(Y), \rH_*(Z)$,
 as de Rham cohomologies of differential forms.
 (It is possible to consider these homologies and cohomologies with
 coefficients in $\Z$ \cite{Fr1}.)
 Recall that the Poincar\'e duality isomorphism
 $\eta_Y:\rH_*(Y)\to \rH^*(Y)$,
 which maps a $k$-cycle to closed $\dim Y-k$ form.  ($\eta_Y^*=\eta_Y$.)
 Then one defines
 $f^*: \rH^*(Y)\to \rH^*(Y)$ and its dual $f_*: \rH_*(Y)\to \rH_*(Y)$
 as
 $$f^*:=\eta_Y\pi_*\eta_Z^{-1}\tilde f^*, \quad f_*:=\tilde
 f_*\eta_Z^{-1}\pi^*\eta_Y.$$
 It can be shown that $f_*,f^*$ do not depend on the resolution of
 $f$, i.e. on $Z$.
 Let $\rho(f_*)=\rho(f^*)$ be the spectral radii
 of $f_*,f^*$ respectively.  (As noted above $f_*,f^*$ can be
 represented by matrix with integer entries.  Hence $\rho(f_*)$ is
 an algebraic integer.)
 Then the \emph{dynamical} spectral radius of $f_*$  is defined as
 \begin{equation}\label{rdyndef}
 \rdyn (f_*)= \limsup_{m\to\infty}
 (\rho((f^m)_*))^{\frac{1}{m}}.
 \end{equation}
 (Note that $\rdyn(f_*)$ is a limit of algebraic integers, so it may
 not be an algebraic integer.)

 Assume that
 $f:Y\to Y$ is holomorphic.  Then $f_*, f_*$ are the standard
 linear maps on homology and cohomology of $Y$.  So
 $(f^m)_*=(f_*)^m, (f^m)^*=(f^*)^m$ and $\rdyn(f_*)=\rho(f_*)$.
 It was shown by the author that $\log\rho(f_*)=h(f)$ \cite{Fr1}.
 This equality followed from the observation that
 $h(f)$ is the \emph{volume growth} induced by $f$.
 View $Y$ as a submanifold of $\P^n$, is endowed the induced
 Fubini-Study Riemannian metric and
 with the induced K\"ahler $(1,1)$ closed form $\kappa$.
 Let $V\subseteq Y$ be any irreducible variety of complex dimension $\dim
 V\ge 1$.  Then the volume of $V$
 is  given by the Wirtinger formula $\vol(V)=\frac{1}{(\dim V)!} \int
 _V \kappa^{\dim V} (=\kappa^{\dim V}(V))$.
 Let $L_k\subset \P^n$ be a linear space of codimension $k$.
 ($L_0:=\P^n$.)
 Assume that $L_k$ is in general position.  Then $L_k\cap V$ is a
 variety of dimension $\dim V -k$.  For $k<\dim V$ the variety $L_k\cap V$
 is irreducible.  For $k=\dim V$ the variety $L_k\cap V$
 consists of a fixed number of points, independent of a generic $L_k$,
 which is called the degree of $V$, and denoted by $\deg V$.
 It is well known that $\deg V=\vol (V)$.
 The homology class of $L_k\cap V$, denoted by $[L_k\cap V]$,
 is independent of $L_k$.  Since
 $\vol(L_k\cap V)$ can be expressed in terms of the cup product
 $\an{[L_k\cap V], [\kappa^{\dim V-k}]}$, or equivalently as $\deg L_k\cap
 V$, this volume is an \emph{integer}, which is independent of
 the choice of a generic $L_k$.  Thus the $j$-th volume growth, of the
 subvariety $L_{\dim Y-j}\cap Y$ of dimension $j$,
 induced by $f$ is given by
 \begin{eqnarray}
 &&\beta_j:=\limsup_{m\to\infty} \frac{\log \an{(f^m)_*[L_{\dim Y-j}\cap Y],
 [\kappa^j]}}{m}, \quad j=1,\ldots,\dim Y, \nonumber \\
 &&H(f):=\max_{j=1,\ldots,\dim Y}
 \beta_j. \label{jvolgrth}
 \end{eqnarray}
 (See \cite[(2)]{Fr1} and \cite[(2.8)]{Fr4}.)  From the well known equality
 $\rho(f_*)=\lim_{m\to\infty}
 ||f_*^m||^{\frac{1}{m}}$, for any norm $||\cdot||$ on $\rH_*(Y)$, it
 follows that $H(f)\le \rho(f_*)$.  Newhouse's result \cite{New} claims that
 $h(f)\le H(f)$.  Combining this inequality with Yomdin's inequality
 \cite{Yom} $h(f)\ge \log\rho(f_*)$ we deduced in \cite{Fr1}:
 \begin{equation}\label{eqenspecr}
 H(f)=\log\rho(f)=h(f),
 \end{equation}
 which is a logarithm of an algebraic integer.

 Let $\rK\subset \rH_{*}(Y)$ be the
 cone generated by the homology classes $[V]$ corresponding to all
 irreducible projective varieties $V\subseteq Y$.
 Note that $f_*(\rK)\subseteq \rK$.  Let
 $\rH_{*,a}(Y):=\rK-\rK \subset \rH_*(Y)$ be the subspace generated by
 the homology classes of projective varieties in $Y$.  Then
 $f_{*}:\rH_{*,a}(Y)\to \rH_{*,a}(Y)$ and denote $f_{*,a}:=f_*|\rH_{*,a}(Y)$.
 Using the theory of nonnegative operators on finite dimensional
 cones $\rK$, e.g. \cite{BP}, it follows that
 $H(f)=\log\rho(f_{*,a})$.

 Assume again that $f:Y \dashrightarrow Y$ is rational dominant.
 Then $f_*(K)\subseteq K$ so $f_*:\rH_{*,a}(Y)\to
 \rH_{*,a}(Y)$ and denote $f_{*,a}:=f_*|\rH_{*,a}*(Y)$.
 Hence we can define $H(f)$, the volume growth induced by $f$, as in
 (\ref{jvolgrth}) \cite{Fr1,Fr4}.  Similar quantities were considered in \cite{RS, BV}.
 It is plausible to assume that $H(f)=\log\rdyn(f_*)$
 and we conjecture a more general set of equalities in the next
 section.

 It was shown in \cite{Fr1} that the results on of Friedland-Milnor
 \cite{FM} imply the inequalities
 \begin{equation}\label{majineq}
 (f^{m})_{*,a} \le (f_{*,a})^{m}  \textrm{ for all }
 m\in\N,
 \end{equation}
 for certain polynomial biholomorphisms of $\C^2$, (which are
 birational maps of $\P^2$.)

 It was claimed in \cite[pp. 367]{Fr1} that if (\ref{majineq}) holds then
 the sequence $(\rho((f^m)_{*,a}))^{\frac{1}{m}}$, $m\in\N$ converges.
 (This is probably
 wrong.  One can show that under the assumption (\ref{majineq}) for all rational
 dominant maps $f:Y\dashrightarrow Y$ one has $\rho((f^q)_{*,a})^p \ge
 \rho((f^{pq})_{*,a})$ for any $p,q\in\N$.)
 It was also claimed in \cite[Lemma 3]{Fr1} that (\ref{majineq}) holds in
 general.  Unfortunately this result is false, and a counterexample
 is given in \cite[Remark 1.4]{Gue}.   Note that if $f:Y \to Y$
 holomorphic then equality in (\ref{majineq}) holds.
 Hence all the results of \cite{Fr1} hold for holomorphic maps.

 \section{Entropy of rational maps}
 Let $f:Y \dashrightarrow Y$ be a rational dominating map.  (We will assume
 that $f$ is not holomorphic unless stated otherwise.)  In order
 to define the entropy of $f$ we need to find the largest subset
 $X\subseteq Y\backslash \Sing f$ such that $f:X \to X$.
 Let $X_k$ is the collection of all $x\in Y$ such that $f^j(x)\in Y\backslash \Sing
 f$ for $j=0,1,\ldots,k$.  Then $X_k$ is open and $Z_k:=Y\backslash
 Y_k$ is a strict subvariety of $Y$.  Clearly $X_k \supseteq
 X_{k+1}, Y_k\subseteq Y_{k+1}$ for $k\in\N$.  Then
 $X:=\cap_{k=1}^{\infty} X_k$ is $G_{\delta}$ set.  Let $\kappa$ the the
 closed $(1,1)$-K\"ahler form on $Y$.  Then $\kappa^{\dim Y}$ is a
 canonical volume form on $Y$.  Hence $\kappa^{\dim Y}(X)=\kappa^{\dim
 Y}(Y)$, i.e. $X$ has the full volume.

 Since $Y$ is a compact Riemannian manifold, $Y$ is a compact metric
 space.  Thus we can define the three entropies
 $h(f),h_B(f),h_F(f)$ in \S3.  So
 $$h_B(f)\le h(f)= h_F(f).$$

 Assume that $f:\C^n\to \C^n$ is a polynomial dominating map.
 Then $f$ lifts to a rational dominating map $f:\P^n \dashrightarrow
 P^n$, which may be holomorphic.  Hence $X\supseteq \C^n$.
 Assume that $f$ is a proper map of $\C^n$.  Recall that one point
 compactification of $\C^n$, denoted by $\C^n\cup\{\infty\}$, is
 homeomorphic to the $2n$ sphere $\rS^{2n}$.  Then $f$ lifts to a
 continuous map $ f_s: \rS^{2n}\to \rS^{2n}$.  Thus we can define
 the entropy $h(f_s)$.  It is not hard to show that $h(f_s)\le
 h_F(f)$.

 Let $\orb X\subset Y^{\N}$ be the orbit space of $f$, and let $\cX$
 be its closure.  $\cX$ is closely related to the graph construction
 discussed in \cite{Gro, Fr1,Fr4,Fr2,Fr3} as well as in other papers.
 Denote by $\Gamma(f)\subset Y^2$ the closure of the set
 $\{(x,f(x)), \;x\in Y\backslash \Sing X\}$ in $Y^2$.
 Then $\Gamma(f)$ is an irreducible variety of dimension $\dim Y$ in
 $Y^2$.
 Note that the projection of $\Gamma$ on the first or second factor
 of $Y$ in $Y^2$ is $Y$.
 Without a loss of generality we may assume that $\Gamma(f)$
 is smooth.

 Otherwise let $\pi:Z \to Y$ be a blow up of $Y$ such
 that $f:Y\dashrightarrow Y$ lifts to a holomorphic map $\tilde f:
 Z\to Y$.  Let $\Gamma_1(f):=\{(z,\tilde f(z)):\; z\in Z\}\subset
 Z\times Y$.  Then $\Gamma_1(f)$ is smooth variety of dimension $\dim Y$.
 Note that $\hat\pi:Z^2\to Z\times
 Y$ given by $(z,w)\mapsto (z,\pi(w))$ is a blow up of $Z\times Y$.
 Lift $\tilde f$
 to $\hat f: Z\dashrightarrow Z$.
 Then $\Gamma( \hat f)\subset Z^2$ is a blow up $\Gamma_1(f)$, hence
 $\Gamma(\hat f)$ is smooth.

 Let $\Gamma\subset Y^2$ be a closed irreducible smooth variety of
 dimension $\dim Y$ such that the projection of $\Gamma$ on the
 first or second component is $Y$.  Define

 \begin{eqnarray*}
 &&Y^k(\Gamma):=\{(x_1,\ldots,x_{k})\in Y^k,\;(x_i,x_{i+1})\in \Gamma \textrm{ for
 } i=1,\ldots,k-1\}, k=2,\ldots,\\
 &&Y^{\N}(\Gamma):=\{(x_1,\ldots,x_{k},\ldots)\in Y^{\N},\;(x_i,x_{i+1})\in \Gamma \textrm{ for
 } i\in\N\}.
 \end{eqnarray*}
 Note that $Y^k(\Gamma)$ is an irreducible variety of dimension
 $\dim Y$ in $Y^k$ for $k=2,\ldots$.
 Note that $Y^{\N}(\Gamma)$ is a $\sigma$ invariant compact subset
 of $Y^{\N}$, i.e. $\sigma(Y^{\N}(\Gamma))\subseteq Y^{\N}(\Gamma)$.
 Let $h(\Gamma)=h(\sigma|Y^{\N}(\Gamma))$.
 $Y$, viewed as a submanifold of $\P^n$, is endowed the induced
 Fubini-Study Riemannian metric and
 with the K\"ahler $(1,1)$ form $\kappa$.
 Then $Y^k$ has the corresponding induced product Riemannian metric,
 and $Y^k$ is K\"ahler, with the $(1,1)$ form $\kappa_k$ .
 Let $\vol(Y^k(\Gamma))=\kappa_k^{\dim Y}(Y^k(\Gamma))$ be volume of the variety
 $Y^k(\Gamma)$.  Then the volume growth of $\Gamma$ is given by
 \begin{equation}\label{lov}
 \lov (\Gamma):= \limsup_{k\to\infty} \frac{\log \vol
 (Y^k(\Gamma))}{k}.
 \end{equation}
 The fundamental inequality due to Gromov \cite{Gro}
 \begin{equation}\label{gromin}
 h(\Gamma)\le \lov(\Gamma).
 \end{equation}
 Since the paper of Gromov was not available to the general public
 until the appearance of \cite{Gro},
 the author reproduced Gromov's proof of (\ref{gromin}) in
 \cite{Fr4, Fr2}.  Using the above inequality Gromov showed that
 $h(f)\le \log\deg f$ for any holomorphic $f:\P^n \to \P^n$.

 Let $f:Y\dashrightarrow Y$ be rational dominating.
 Then $\cX=Y^{\N}(\Gamma(f))$.  Hence
 \begin{equation}\label{frenteq}
 h_F(f)=h(\Gamma(f)).
 \end{equation}
 If $\Gamma(f)$ is smooth then Gromov's inequality yields that
 \begin{equation}\label{grominf}
 h_F(f)\le \lov(f):=\lov(\Gamma(f)).
 \end{equation}

 \begin{con}\label{eqalentr}  Let $Y$ be a smooth projective variety
 and $f:Y\dashrightarrow Y$ be a rational dominating map.
 Then there exists a smooth projective variety $\hat Y$ and a
 birational map $\iota:Y \dashrightarrow \hat Y$, such that the lifting $\hat f:
 \hat Y \dashrightarrow \hat Y$ satisfies
 (\ref{entropconeq}).
 \end{con}

 We now review briefly certain notions, results and conjectures in
 \cite[S3]{Fr2}.
 Let $\Gamma\subset Y^2$ be as above, and denote by $\pi_i(\Gamma)\to Y$
 the projection of $\Gamma$ on the
 $i$-th component of $Y$ in $Y\times Y$ for $i=1,2$.
 Since $\dim\Gamma=\dim Y$ and $\pi_1(\Gamma)=\pi_2(\Gamma)=Y$,
 then $\deg \pi_i$ is finite and $\pi_i^{-1}(y)$ consists of exactly
 $\deg\pi_i$ distinct points for a generic $y\in Y$ for $i=1,2$.
 One can define a
 linear map $\Gamma_{*}:\rH_*(Y)\to \rH_{*}(Y)$ given by
 $\Gamma_*: \pi_1^*\eta_{\Gamma}^{-1}\pi_2^*\eta_Y$.
 (This is an analogous definition of $f_*$, where
 $f:Y\dashrightarrow Y$ is dominating.)  One can show that
 $\Gamma_*(\rH_{*,a}(Y))\subseteq \rH_{*,a}(Y)$.  Let
 $\Gamma_{*,a}:=\Gamma_*|\rH_{*,a}(Y)$.

 $\Gamma\subset Y^2$ is called a \emph{proper} if each $\pi_i$ is
 finite to one.
 Assume that $\Gamma$ is proper.
 Then
 \begin{equation}\label{gamineq}
 \log\rho(\Gamma_{*,a})\ge
 \lov(\Gamma).
 \end{equation}
 It is conjectured that for a proper $\Gamma$
 \begin{equation}\label{gamcon}
 \log\rho(\Gamma_{*,a})=
 \lov(\Gamma)=h(\Gamma).
 \end{equation}
 Note that if $f:Y\to Y$ is dominating and holomorphic then
 $\Gamma(f )$ is proper, $\Gamma_{*,a}=f_*|\rH_{*,a}(Y)$
 and the above conjecture holds.

 We close this section with observations and remarks which are not
 in \cite{Fr2}.
 Assume that $f:Y\dashrightarrow Y$ be a rational dominating and
 $Z:=\Gamma(f)\subset Y^2$ smooth.  Then $\pi_1:\Gamma(f)\to Y$ is a
 blow up of $Y$, and $\pi_2:\Gamma\to Y$ can be identified with
 $\tilde f: Z\to Y$.  It is straightforward to show that
 $f_*=\Gamma( f)_*$.

 It seems to the author that the arguments given in \cite[Proof Thm
 3.5]{Fr2} imply that (\ref{gamineq}) holds for any smooth variety
 $\Gamma\subset Y^2$ of dimension $\dim Y$ such that
 $\pi_1(Y)=\pi_2(Y)=Y$.  Suppose that this result is true.
 Let $f:Y \dashrightarrow Y$ be rational and dominating.  Assume
 that $\Gamma(f)\subset Y^2$ is smooth.  Then (\ref{gamineq}) would
 imply that $\lov(f) \le \log\rho(f_*)$.  Applying the same
 inequality to $(f^k)_*$ and combining it with (\ref{grominf}) one would able to deduce:
 \begin{equation}\label{lovrdynin}
 h_F(f)\le \lov(f)\le \log\rdyn (f_*).
 \end{equation}

 \section{Currents}

 Many recent advances in complex dynamics in several complex variables were
 achieved using the notion of a \emph{current}.  See for example the
 survey article \cite{Sib}.  Recall that on an $m$-dimensional
 manifold $M$ a current of degree $m-p\ge 0$ is a linear functional
 on all smooth $p$-differential forms $\cD^p(M)$ with a compact
 support, where $p$ is a nonnegative integer.

 Let $f:Y\dashrightarrow Y$ be a meromorphic dominating self map of
 a compact K\"ahler manifold of complex dimension $\dim Y$, with the $(1,1)$ K\"ahler
 form $\kappa$.
 Let $f^*\kappa$ be a pullback of $\kappa$.  Then
 $f^*\kappa$ is a current on $Y\backslash \Sing f$.
 Define the \emph{p-dynamic degree} of $f$ by
 $$\lambda_p(f):=\limsup_{k\to\infty}\big(\int_{Y\backslash \Sing f^k} (f^k)^*\kappa^p\wedge
 \kappa^{\dim
 Y-p}\big)^{\frac{1}{k}},\quad p=1,\ldots,\dim Y.$$
 It is shown in \cite{DS} that
 the dynamical degrees are invariant with respect to a bimeromorphic
 map $\iota:Y\dashrightarrow Z$, where $Z$ is a compact K\"ahler
 manifold. (See also \cite{Gue} for the case where $Y,Z$ are projective varieties.)
 Moreover
 \begin{equation}\label{DSineq}
 \lov(f)\le \max_{p=1,\ldots,\dim Y}
 \log \lambda_p(f).
 \end{equation}
 Assume that $Y$ is
 a projective variety.  It can be shown that the dynamic degree $\lambda_p(f)$
 is equal to $e^{\beta_{\dim Y-p}}$ for $p=1,\ldots,\dim Y$, which are defined in
 (\ref{jvolgrth}), where $\beta_0:=\beta_{\dim Y}$.
 Hence
 \begin{equation}\label{Hfchar}
 H(f)=\max_{p=1,\ldots,\dim Y}
 \log \lambda_p(f),
 \end{equation}
 where $H(f)$ is defined in (\ref{jvolgrth}).
 Thus $H(f)$ can be viewed as the \emph{algebraic entropy}
 of $f$ \cite{BV}.  \cite[Lemma 4.3]{GS}  computes $H(f)$
 for a large class of automorphisms of $\C^k$, and see also \cite{DS2, Gue}.
 Combine (\ref{grominf}) with (\ref{DSineq}) and (\ref{Hfchar}) to
 deduce the inequality $h_F(f)\le H(f)$, which was conjectured in
 \cite[Conjecture 2.9]{Fr4}.

 Consider the following example
 $f:\C^2\to \C^2, (z,w)\mapsto (z^2,w+1)$ \cite[Example 1.4]{Gue2}.
 Since $f$ is proper we have $f_s:S^4\to S^4$.  Clearly $S^4$ is the
 domain of attraction of the fixed point $f_s(\infty)=\infty$.
 Hence $h(f_s)=0$.  Lift $f$ to $f:\P^2\dashrightarrow \P^2$.  Then
 $f$ has a singular point $\a:=(0,1,0)$ and any other point at the line
 at infinity $(1,w,0)$ is mapped to a fixed point $\b:=(1,0,0)$.
 So $X=\P^2\backslash\{a\}$, and
 $\Gamma(f)=\{(\z,f(\z)):
 \z\in\P^2\backslash\{a\}\}\cup\{(\a,(z:w:0)):(z:w)\in\P\}$, which is equal to the blow up
 of $\P^2$ at $\a$.  On $(\P^2)^{\N}(\Gamma(f))$ $\sigma$ has two
 fixed points: $(\a,\a,\ldots), (\b,\b,\ldots)$.  The set
 $\cX_0:=((\x,f(\x),\ldots,): \x\in A_0:=\{(z,w,1), |z|\le 1\}\}$
 is in the domain of the attraction of $(\a,\a,\ldots)$.
 The set $(\P^2)^{\N}(\Gamma(f))\backslash
 (\cX_0\cup\{(\a,\a,\ldots)\}$ is in the domain of the attraction of
 $(\b,\b,\ldots)$.  Hence $h(f)=0$.
 Observe that $\hat f:(\P\times\P)\to (\P\times \P)$, given as $((z:s),(w:t))\mapsto
 ((z^2:s^2),(w+t:t))$, is the lift of $f$ to $(\P\times\P)$.  $\hat
 f$ is holomorphic and $h(\hat f)=H(\hat f)=\log 2$.  Since $\P\times \P$ is birational
 to $\P^2$ it follows that $H(f)=\log 2 > h_F(f)=0$.  In particular
 $h_F(f)$ is not a birational invariant \cite{Gue2}.
 Note that
 Conjecture \ref{eqalentr} is valid for this example.
 Additional examples in \cite{Gue1, Gue2} support the Conjecture
 \ref{eqalentr}.

 Assume now that $f:\C^2 \to \C^2$ is a polynomial
 automorphism, hence $f$ is proper.  It is shown in \cite{FM} that $h(f_s)=h(f,K)$
 for some compact subset of $\C^2$.  Furthermore the results of
 \cite{FM} and \cite{Smi} imply that $h(f,K)=H(f)$.  One easily
 deduce that $H(f)=\rdyn(f_*)$.  Clearly $h_B(f)\ge h(f,K)$.
 Then the inequalities $h_F(f)\le \lov(f)\le H(f)$ yield Conjecture
 \ref{eqalentr}.  See \cite{BD,DS1,HP} for additional results on
 entropy of certain rational maps.

 The inequality (\ref{DSineq}) and
 its suggested variant (\ref{lovrdynin}) can be viewed as Newhouse
 type upper bounds \cite{New} which shows that the volume growth
 bounds from above the entropy of a rational dominating map.
 In order to prove Conjecture \ref{eqalentr} one needs to prove
 a suitabe Yomdin type lower bound  for the entropy of $f$.


\begin{thebibliography}{99}
 \bibitem{AKM} R.L. Adler, A.G. Konheim and M.H. McAndrew,
 Topological entropy, \emph{Trans. Amer. Math. Soc.} 114 (1965),
 309--311.
 \bibitem{BD} E. Bedford and J. Diller, Real and complex dynamics of a
 family of birational maps of the plane: the golden mean subshift,
 \emph{Amer. J. Math.}  127 (2005), 595--646.
 \bibitem{BV} M.P. Bellon and C.-M. Viallet,  Algebraic entropy,
 \emph{Comm. Math. Phys.}  204 (1999), 425--437.
 \bibitem{BP} A. Berman and R.J. Plemmons, \emph{Nonnegative
 Matrices in the Mathematical Sciences}, Academic Press, 1979.
 \bibitem{Buf} A. Bufetov, Topological entropy of free semigroup actions and skew-product
 transformations, \emph{J. Dynam. Control Systems}  5 (1999), 137--143.
 \bibitem{Can} S. Cantat, Dynamique des automorphismes des surfaces
 $K3$, Acta Math.  187  (2001), 1--57.
 \bibitem{Din}  E.I. Dinaburg, A correlation between topological entropy and metric
 entropy, \emph{Dokl. Akad. Nauk SSSR} 190 (1970), 19--22.
 \bibitem{DS} T.-C. Dinh and N. Sibony, Regularization of currents and
 entropy, \emph{Ann. Sci. École Norm. Sup.} 37 (2004), 959--971.
 \bibitem{DS1} T.C. Dinh and N. Sibony, Green currents for holomorphic automorphisms of compact
 K\"ahler manifolds, \emph{J. Amer. Math. Soc.} 18 (2005), 291--312.
 \bibitem{DS2} T.-C. Dinh and N. Sibony, Une borne supérieure pour
 l'entropie topologique d'une application rationnelle,
 \emph{Ann. of Math.} 161 (2005), 1637--1644.
 \bibitem{Fis}  A. Yu. Fishkin,
 An analogue of the Misiurewicz-Przytycki theorem for some mappings,
 \emph{Uspekhi Mat. Nauk 56} 337 (2001), 183--184.
 \bibitem{Fr1}  S. Friedland, Entropy of polynomial and rational
 maps, \emph{Ann. of Math.} 133  (1991), 359--368.
 \bibitem{Fr4} S. Friedland, Entropy of rational self-maps of projective varieties,
 \emph{International Conference on Dynamical Systems and Related Topics},
editor: K. Shiraiwa, Advanced Series in Dynamical Systems, vol. 9,
128-140, World Scientific Publishing Co., Singapore 1991.
 \bibitem{Fr2} S. Friedland, Entropy of algebraic maps, Proceedings of the
 Conference in Honor of Jean-Pierre Kahane, \emph{J. Fourier Anal. Appl.}  1995,
 Special Issue, 215--228.
 \bibitem{Fr3} S. Friedland, Entropy of graphs, semigroups and groups, Ergodic theory of $Z\sp d$
 actions (Warwick, 1993--1994),  319--343,
 \emph{London Math. Soc. Lecture Note Ser.} 228, Cambridge Univ. Press, Cambridge, 1996.
 \bibitem{FM} S. Friedland and J. Milnor, Dynamical properties of plane polynomial automorphisms,
 \emph{J. Ergod. Th. \& Dynam. Sys.} 9 (1989), 67--99.
 \bibitem{Gdm} T.N.T. Goodman, Relating topological entropy and measure entropy.
 Bull. London Math. Soc.  3 (1971), 176--180.
 \bibitem{Goo} L.W. Goodwyn, Topological entropy bounds measure-theoretic entropy.
 \emph{Proc. Amer. Math. Soc.}  23 (1969), 679--688.
 \bibitem{GH} P. Griffiths and J. Harris, \emph{Principles of Algebraic
 Geometry}, Wiley Interscience, 1978.
 \bibitem{Gro} M. Gromov, On the entropy of holomorphic maps,
 \emph{Enseign. Math.} 49 (2003),  217--235.
 \bibitem{Gue1} V. Guedj, Courants extr\'emaux et dynamique complexe.
 \emph{Ann. Sci. École Norm. Sup.} 38 (2005), 407--426.
 \bibitem{Gue2} V. Guedj, Entropie topologique des applications
 m\'eromorphes, \emph{Ergodic Theory Dynam. Systems}  25 (2005), 1847--1855.
 \bibitem{Gue} V. Guedj, Ergodic properties of rational mappings
 with large topological degree, \emph{Ann. of Math.} 161  (2005), 1589--1607.
 \bibitem{GS} V. Guedj and N. Sibony, Dynamics of polynomial automorphisms of
 $\bold C\sp k$, \emph{Ark. Mat.}  40 (2002), 207--243.
 \bibitem{HNP} B. Hasselblatt, Z. Nitecki and J. Propp,
 Topological entropy for non-uniformly continuous maps,
 arXiv:math.DS/0511495 v1, 20 Nov. 2005.
 \bibitem{HP} B. Hasselblatt and J. Propp, Monomial maps and algebraic
 entropy, arXiv:mathDS/0604521 v1, 25 Apr. 2006.
 \bibitem{Kol} A.N. Kolmogorov, A new metric invariant of transitive
 dynamical systems and Lebesgue space automorphisms, \emph{Dokl. Acad. Sci. USSR}
 119 (1958), 861--864.
 \bibitem{Lyu}  M.Yu. Lyubich, Entropy of analytic endomorphisms of the Riemann sphere.
 \emph{Funktsional. Anal. i Prilozhen.}  15 (1981), 83--84.
 \bibitem{McM} C.T. McMullen, Dynamics on $K3$ surfaces: Salem numbers and Siegel
 disks, \emph{J. Reine Angew. Math.}  545 (2002), 201--233.
 \bibitem{MP} M. Misiurewicz and F. Przytycki, Topological entropy and degree of smooth mappings.
 \emph{Bull. Acad. Polon. Sci. Sér. Sci. Math. Astronom. Phys.}  25, (1977), 573--574.
 \bibitem{New} S.E. Newhouse, Entropy and volume,  \emph{Ergodic Theory Dynam. Systems } 8$\sp *$
 (1988),  Charles Conley Memorial Issue, 283--299.
 \bibitem{Prz} F. Przytycki, An upper estimation for topological entropy of diffeomorphisms.
 \emph{Invent. Math.}  59  (1980), 205--213.
 \bibitem{RS} A. Russakovskii and B. Shiffman,  Value distribution for sequences
 of rational mappings and complex dynamics.
 \emph{Indiana Univ. Math. J.} 46 (1997), 897--932.
 \bibitem{Shu}  M. Shub, Dynamical systems, filtrations and
 entropy, \emph{Bull. Amer. Math. Soc.} 80  (1974), 27--41.
 \bibitem{Sib} N. Sibony, Nessim Dynamique des applications rationnelles de $\bold P\sp
 k$, \emph{Panor. Synthèses} 8 (1999), 97--185, Soc. Math. France, Paris.
 \bibitem{Smi} J. Smillie, The entropy of polynomial diffeomorphisms of
 $C\sp 2$, \emph{Ergodic Theory Dynam. Systems} 10 (1990), 823--827.
 \bibitem{Wal} P. Walters, \emph{An Introduction to Ergodic Theory},
 Springer-Verlag, 1982.
 \bibitem{Yom} Y. Yomdin, Volume growth and entropy, \emph{Israel J. Math.}  57 (1987), 285--300.



 \end{thebibliography}
\end{document}